\documentclass[12pt,leqno]{article}
\usepackage{amstext}
\usepackage[dvips]{graphics}
\usepackage{theorem}
\usepackage{amssymb}
\usepackage{latexsym}
\newtheorem{thm}{Theorem.}[section]
\newtheorem{prop}[thm]{Proposition.}
\newtheorem{cor}[thm]{Corollary.}
\newtheorem{lem}[thm]{Lemma.}
\newtheorem{rem}[thm]{Remark.}

\newtheorem{defn}[thm]{Definition.}

\newtheorem{prop-def}[thm]{Proposition-Definition.}

\def\pmb#1{\setbox0=\hbox{#1}%
\kern-.025em\copy0\kern-\wd0
\kern.05em\copy0\kern-\wd0
\kern-.025em\raise.0433em\box0}

\newcommand{\lra}{\longrightarrow}

\newcommand{\be}{\begin{enumerate}}
\newcommand{\ee}{\end{enumerate}}

\newcommand{\br}{\begin{array}}
\newcommand{\er}{\end{array}}

\newcommand{\vg}{>\!\!>}
\newcommand{\up}{^{\prime}}
\newcommand{\upp}{^{\prime\prime}}

\newcommand{\frap}[2]{\mbox{$\frac{#1}{#2}$}}

\newcommand{\binom}[2]{\mbox{$\left(\!\br{c}{#1}\vspace{-.5ex}
\\{#2}\er\!\right)$}}

\newfont{\ninemsbm}{msbm10 scaled 0900}
\newfont{\tenmsbm}{msbm10 scaled 1100}
\newfont{\nineeufb}{eufb10 scaled 0900}
\newfont{\teneufb}{eufb10 scaled 1100}
\newfont{\teneusm}{eusm10 scaled 1200}
\newfont{\nineeusm}{eusm10 scaled 0900}
\newfont{\bal}{cmmib10 scaled 0900}
\newfont{\tencmmib}{cmmib10 scaled 1000}

\newcommand{\gotm}[1]{\mbox{{\teneufb {#1}}}}
\newcommand{\sbm}[1]{\mbox{{\tenmsbm {#1}}}}

\newcommand{\fasm}[1]{\mbox{{\teneusm {#1}}}}

\newcommand{\bm}[1]{\mbox{{\boldmath ${#1}$}}}

\begin{document}

\begin{center}{\huge\bf An asymptotic vanishing
theorem for generic unions of multiple points
\vspace{3ex}}\\
J.Alexander\hspace{3cm}A. Hirschowitz\end{center}
\section{Introduction}

This work is devoted to the following  asymptotic statement :
\begin{thm}\label{thmm} Let $X$ be a projective geometrically reduced
and irreducible scheme over a field $k$ of (arbitrary) characteristic
$p$ and let $\fasm M\, ,\, \fasm L$ be line bundles on $X$ with $\fasm
L$ ample. If $p$ is positive then suppose further that $X$ is smooth in
codimension one. For fixed $m\geq 0$ there exists $d_0=d_0(m)$,
depending also on $X,\fasm L,\fasm M$, such that for any $d\geq d_0$
and any generic union $Z$ of (fat) points of multiplicity $\leq m$ the
canonical map
$$H^0(X,\fasm M\otimes\fasm L^d)\lra H^0(Z,\fasm O_Z
\otimes\fasm M\otimes \fasm L^d)$$
has maximal rank.\end{thm}

Here, as usual, we call (fat) point of multiplicity $m$ in $X$, any
subscheme defined by $\fasm I_z^m$, where $\fasm I_z$ is the ideal
sheaf of a point $z$ in the smooth locus of $X$. The reader may prefer
the following statement,  which is more or less equivalent to the
preceding one:

\begin {cor}Let $X$, $\fasm M\, ,\, \fasm L$, $m$ be as above. There
exists an integer $\ell$
such that for any generic union $Z$ of (fat) points of multiplicity at
most $ m$ and  of  total degree (i.e length)
at least $\ell$, all the canonical maps
$$H^0(X,\fasm M\otimes\fasm L^d)\lra H^0(Z,\fasm O_Z\otimes\fasm
M\otimes \fasm L^d)$$
have maximal rank. \end {cor}
Note that the above statement applies as soon as the
number of points is at least $\ell$.

To simplify the presentation and highlight the essential elements, a
detailed proof will only be given in the case $\fasm M=\fasm O_X$.
The easy modifications needed to prove the general result are then
outlined in \S 7 along with another variant.
\begin{rem} The statement of the theorem is false for $p>0$ if we allow
$X$ to be singular in codimension one. This is illustrated in the
example \ref{cp}.
\end{rem}

These results are already new in the case where $X$
is the projective plane (with $\fasm M = \fasm O$ and $ \fasm
L=\fasm O (1)$). Indeed, even in that case, the expected
vanishing theorem for generic unions of fat points [S, Ha, Hi2] is
still unproven, see a survey in [G] and more recent contributions in
[Xu, ShT, CM1, CM2, M]. Reformulations of the general problem and its
relation to other topics have been considered at length in [N, I1,
I2, MP]. Much attention has been paid to the ``homogeneous'' case on
$\sbm P^n$, namely when all the points have the same multiplicity $m$:
see  [AC] or [Hi1] for $n=m=2$, [A, AHi1,2,3] for $m=2$ and $n$
arbitrary,  [Hi1] for $m=3, n=2, 3$, [LL2] for  $m=4, n=2$ in
perfectly adjusted cases, and finally [CM1, CM2], where they have settled
completely the equal multiplicity cases $m\leq 12, n=2$ by a new and promising
method. Concerning the heterogeneous plane case, we can just
mention the recent work of Th. Mignon [M], where the case
of multiplicities at most four is
completely elucidated, using
our differential Horace lemma presented below.

In [AHi2] we developed a new technique of a differential nature for
the case $m=2$ which, in that and later papers,  made it possible not
only to  solve some delicate low-degree cases in [AHi2,3], but
also to simplify the proof for the high-degree case [AHi4][C]. The main
new ingredient in the proofs of the present paper is  an extension of
this technique applicable to higher-order fat points ($m> 2$), see
lemma \ref{diff} and \ref{prin.cor}. The new lemma does not
imply the
multiplicity two lemma of [AHi2], and an entirely new proof is needed.

In sections 2-7, we present the proof of the theorem. Sections 8-9 are
devoted to our differential Horace lemmas. Indeed, the results
presented there (see \S 9) are substantially more general than
\ref{diff}. While for the present asymptotic statement, \ref{diff} is
perfectly sufficient, the full strength  of \ref{prin.cor} will be
much more efficient for concrete cases with small $n$ and $m$. The
proof of
\ref{prin.cor} is achieved by an ideal theoretic argument. We would
like to point out that our original proof of the lemmas computed the
first non-zero derivative of a determinant in a way which owed much to
[LL1,2].

\noindent{\bf Outline of the proof of the theorem}

In the remainder of this introduction we will try to illustrate
the general ideas in the proof of the main theorem in the particular
case of the projective plane. We start with a given maximum
multiplicity $m$ and a sufficiently large degree $d$. We want to prove
a maximal rank statement for a generic union $Z$ of multiple points,
which, by adding simple points we can suppose to be of total degree at
least $(d+2)(d+1)/2$. Horace's method amounts to specialising some of
these points to the generic curve $\Gamma$ of some intermediate
degree $\gamma$. Modulo an analogous maximal rank statement on
$\Gamma$ which we suppose to hold inductively,
our differential lemma can then  be
applied under certain numerical conditions (holding for large $d$) and
we reduce to a new subscheme $\gotm D^{(1)}(Z)$ (the {\it derivative}
of $Z$, see \S 4) and a new degree $d-\gamma$. This  can be safely
applied as long as the current degree, $d_c$, is not too small, say $d_c >
\overline{d}$. But when $d_c$ becomes smaller than or equal to $
\overline{d}$, we have to backtrack in order to complete the proof.
Our trick is to modify the procedure early on so as to generate in the
current subscheme $Z_c$ a sufficient number of unconstrained points of
multiplicity at most $m-1$ (of total degree at least
$(\overline{d}+2)(\overline{d}+1)/2$). So that when the degree of the
current subscheme has been lowered under $\overline{d}$, only points of
multiplicity at most $m-1$ remain. Having chosen $\overline{d}$ large
enough (i.e. $\overline{d}\geq d_0(m-1)$ in the notation of the
theorem), we conclude by induction on $m$.

It remains to explain how we generate these free points (see \S 6): our
differential lemma generates in  $\gotm D^{(1)}(Z)$ a lot of points of
multiplicity smaller than $m$, but all of them lie on  the exploited
divisor $\Gamma$ of degree $\gamma$. The trick here consists in
specialising $\Gamma$ to the union of two generic divisors $\Gamma\up$
and $\Gamma\upp$ of degrees $\gamma\up$ and $\gamma\upp$, with  the
desired number of points specialised to say $\Gamma\up$. If this
number of points is sufficiently small with respect to $\gamma\up$,
these points suffer no constraint by being supported on a curve of
degree $\gamma\up$ and are thus freed. Of course, the points remaining
on $\Gamma\upp$ should not be too numerous, and we have to find
numbers  $d_0, \gamma\up$ and $\gamma\upp$ satisfying all the
necessary inequalities.

A slight complication arises with the
degree of the  current divisor $\Gamma_c$. Indeed,
the number of free points to be generated is computed in terms of the
degree of the divisor which appears at the final stage of the procedure
(this degree must be sufficiently large to comply with the induction
hypothesis). On the other hand, the initial degree $\gamma$ of the
current divisor must be significantly  larger to allow the production
of enough free points. This compels us to lower the degree of the
current divisor, by specialization, at each stage of the procedure
(see \S 5).
\section{The simplified differential lemma}
Throughout this section, $X$ stands for a quasi-projective variety
which is geometrically reduced and irreducible, of dimension $n+1$
over a field $k$ of arbitrary characteristic. Since all statements are
``generic'' one can safely suppose $k$ algebraically closed. The
hypothesis `$X$ is smooth in codimension one if $char(k)>0$' will not
come into play until the proof of the theorem in \S 7.

In this section, we present a weakened form of the differential
lemma which we prove in \S 9, this form being sufficient for our main
theorem. As we already outlined in the previous section, the theorem
is proved by a Horace induction argument. In such an argument,
specialisation techniques are used to place a certain number of points
on a chosen divisor $H$, then the induction hypotheses are applied to
the trace and the residual as defined in the
\begin{defn}\label{trace-residual}
Let $H$ be a Cartier divisor on $X$ and let $W$ be a closed subscheme
of $X$.

The schematic intersection
$$ W^{\prime\prime} =H\cap W $$
defined by the ideal $\fasm I_{H,W^{\prime\prime}}=(\fasm
I_H+\fasm I_W) /\fasm I_H$ of $\fasm O_H$ is called the trace of $W$
on $H$ and  denoted by $\mbox{Tr}_H(W)$ or simply $W^{\prime\prime}$ if
no  confusion is possible.

The closed subscheme of $X$ defined by the conductor ideal
$\fasm I_{W^{\prime}}=(I_W :I_H)$ is called the residual of $W$ with
respect to $H$ and denoted by $\mbox{Res}_H(W)$ or
$W^{\prime}$.

The canonical exact sequence
\begin{equation}\label{res.exact}0\lra \fasm I_{W^{\prime}}(-H)\lra
\fasm I_W
\lra
\fasm I_{H,W^{\prime\prime}}\lra 0\end{equation}
is called the residual exact sequence of $W$ with respect to $H$.
\end{defn}
%

\subsection{Geometric intuition for the differential lemma}

Here we try to share with the reader our intuition for our
differential lemma. Suppose that $X$ is projective and let $\fasm L$ be
a line bundle on $X$. We will keep the notation of the definition in
the remainder of the discussion.

The first thing one needs to take note of is that any basic Horace
type argument is based on the following trivial consequence of the
residual exact sequence (\ref{res.exact}):
$$ \mbox {if $h^0(X,\fasm
I_{H,W\upp}\otimes \fasm L|_H)=0$ and
$h^0(X,\fasm I_{W\up}\otimes
\fasm L(-H))=0$
then $h^0(X,\fasm I_W\otimes \fasm L)=0$}$$
For this to apply, one must have aprori
$\mbox{deg}W\upp \ge h^{0}(H,\fasm L_H)$ and
$\mbox{deg}W\up \ge h^{0}(X,\fasm L(-H))$. In fact to be generally
applicable in
an induction argument, the stronger requirement $\mbox{deg}W\upp =
h^{0}(H,\fasm
L_H)$ is needed. We will therefore say that $h^{0}(H,\fasm
L_H)$ is the {\em critical degree}.

In practice one starts with some general union $G$ of multiple points, then, by
specialising them one by one to the chosen divisor $H$ one hopes to obtain a
specialisation $W$ for which the trace has the critical degree. Since each
point
specialised to $H$ increases the degree of the trace by at least the
multiplicity
of the point, it is not always possible to get exactly the critical degree
using
this process. This is the technical obstacle that the differential lemmas
\ref{prin.lem} and \ref{diff} are designed to overcome.

To see how this comes about, it is enough to consider that $H$ is a line in the
affine plane $X$. The ideal of a point $Z$ of multiplicity $r$ at the origin is
then
\begin{equation}\label{ideal}\gotm (x,y)^r =\gotm n^r \oplus\gotm n^{r-1} y
\oplus
\cdots \oplus
\gotm n y^{r-1}\oplus (y^r)\end{equation}
where $\gotm n$ is the ideal $(x)\subset k[x]$, and each $\gotm n^i$ is the
ideal
of a point of multiplicity $i$ in $H$. In particular the trace corresponds to
$\gotm n^r$. One can then consider that $Z$ is formed by infinitesimally
piling up
the subschemes of $H$ with ideals $\gotm n^i$. We then refer to these
subschemes
of $H$ as the layers of $Z$. Of course only the trace given by $\gotm n^r$ is
actually contained in $H$, the others only appear in successive infinitesimal
neighbourhoods of $H$. Now if we consider $Z$ as the limit of a multiple point
that is translated to the origin along the $y$-axis, it's the layer of highest
multiplicity that arrives in $H$ (or, as might be said, arrives first) and the
degree of the trace increases by $r$.

In the differential approach, one or more points are translated to as many
distinct points supported in $H$. The rate of approach may differ, but all
arrive at the same time. Our corollary \ref{prin.lem} says that if some
sequence of layers, one from each point, have degrees adding up to the
critical degree, then one can consider that these arrive first and
then take their union as the (differential) trace , while the
subsequent remainder becomes the (differential) residual. Precisely,
with respect to the ideal (\ref{ideal}), if the layer corresponding to
$\gotm n^p$ is taken at that point to be its contribution to the
(differential) trace, then the (differential) residual at that point
is the subscheme of the plane defined by the ideal
$$\gotm n^r\oplus\gotm n^{r-1}y\oplus \cdots \oplus\gotm
n^{p+1}y^{r-p-1}\oplus\gotm n^{p-1}y^{r-p}\oplus\cdots
\gotm n y^{r-2}\oplus (y^{r-1}).$$
obtained by {\em slicing off} the corresponding layer. If the
cohomology vanishes as before when the trace and residual are replaced
by the chosen differential versions, then the lemma says that the
cohomology vanishes for $\fasm I_G\otimes \fasm L$. The conclusion now
concerns the general union $G$ and not the specialisation
$W$.

\begin{center}\includegraphics{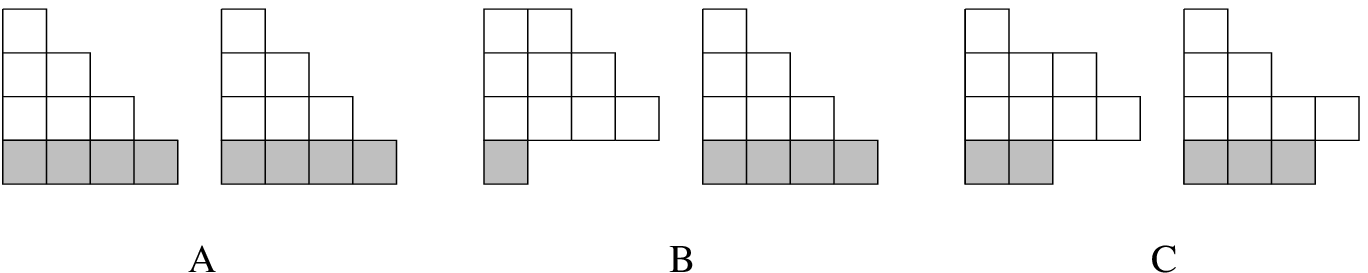}\end{center}

Figure 1. illustrates an example where $X$ is the affine
plane and $H$ is a line, while the critical degree is supposed to be
five. Example A shows two points of multiplicity four in the plane.
The shaded region represents the trace while the unshaded region
represents the residual with respect to the line $H$. From the
standard specialisation point of view, these points are translated, one
by one, to $H$ giving a trace of degree four,  then eight, so that five
is unattainable. The examples B and C show two possibilities for
choosing the differential traces (shaded part) so that the critical
degree is obtained. The differential residuals correspond to the unshaded
part.
\subsection{The simplified lemma}

The simplication of the following lemma with respect to that in
\S 9, is as follows :  in the process just described, instead of
choosing arbitrary slices, we will systematically take the smallest
non-trivial one, which is just a simple point of $H$. In this case, the
(differential) residual falls within the bounds of the following
definition.
\begin{defn}\label{simple residue} Let $H$ be a reduced Cartier
divisor on $X$ and let $z$ be a non-singular point of $H$. We define
the $m^{th}$ simple residue; denoted $D_{H,m}(z)$ or
$D_{m}(z)$ if no confusion can arise; to be the trace of $z^m$
on $(m-1)H$;
$$D_m(z)=z^m\cap H^{m-1}.$$
We will say that $m$ is the {\em multiplicity} of the simple residue.
\end{defn}
With this definition, our simplified lemma, which will be proved in \S9,
reads as follows.
\begin{lem}\label{diff} Suppose $X$ is  projective and furnished with a
line bundle $\fasm L$, and let $H$ be a reduced and irreducible
effective  Cartier divisor on $X$. Let $Z_0$ be a zero-dimensional
subscheme of $X$, and let $a$, $d$ be positive integers. We suppose
that
$$r\; =\; h^0\left(H,\fasm
L|_H\right) - \mbox{deg}\left(Tr_{H}(Z_0)\right) \geq 0$$
and that $m_1,\ldots ,m_r$ are positive integers satisfying
$$\mbox{deg}(Z_0)+\sum_{i=1}^r\,\mbox{$\binom{m_i+n}{n+1}$}\;
\geq \; h^0(X,\fasm L)$$
Let $P_1,\ldots ,P_r$ be generic points in $X$ and $Q_1,\ldots ,Q_r$
be generic points in $H$.
In the notation of \ref{trace-residual} and \ref{simple residue}, set
$$T=Z_0\cup P_1^{m_1}\cup \cdots \cup P_r^{m_r}\quad ;\quad
T\up_{\star}=Z_0\up\cup D_{m_1}(Q_1)\cup \cdots \cup
D_{m_r}(Q_r) \quad ; \quad T\upp_{\star} = Z_0\upp\cup Q_1\cup
\cdots \cup Q_r$$
Then $H^0(X , \fasm I_{T}\otimes\fasm L)=0$ holds as soon as the following two
conditions are satisfied :\vspace{2ex}\\
\makebox[5cm][l]{(dime)}\mbox{$\displaystyle\br{lll}H^0(H,\fasm
I_{T\upp_{\star}}\otimes \fasm L|_H)&=&0\er$}\vspace{2ex}\\
\makebox[5cm][l]{(degue)}\mbox{$\displaystyle \br{lll} H^0(X,\fasm
I_{T\up_{\star}}\otimes \fasm L(-H))&=&0.\er$}
\end{lem}
\begin{rem} The dime and degue concern respectively the differential
trace and the differential residual as discussed above.
\end{rem}
%
\section{\label{CC}Configurations and candidates}
Here we introduce the general class of subschemes of $X$ which we will be
dealing with. From here on, $X$ is projective of dimension $n+1$ and
furnished with an ample line bundle $\fasm O(1)$ of degree $\nu$. We
let $\alpha_0$ be the least integer such that $\fasm O(a)$ is very
ample for $a\geq \alpha_0$ and it will henceforth be understood that
$a\geq \alpha_0$.
\begin{defn}\label{candidates} Let
$G_a$ be the generic effective divisor in
the linear system $\bigl|H^0(X,\fasm O(a))\bigr|$. A $G_a$-residue or
just {\bf residue}, will be any point or simple residue (see
\ref{simple residue}) with support in $G_a$. The {\bf multiplicity}
of a residue will be its multiplicity as a point, or as a simple residue
\ref{simple residue}, respectively.
\newline
Given positive integers $a,m$ an $(a,m)$-{\bf configuration} will be
any subscheme $Z$ of $X$ which is a generic union of points of
multiplicity at most $m$ in $X$, called the free part of $Z$ and
denoted $\mbox{Free}(Z)$, with a generic set of $G_a$-residues equally
of multiplicity at most $m$, called the constrained part of $Z$ and
denoted $\mbox{Const}(Z)$.
\newline
Given a positive integer $d$, we say that an $(a,m)$-configuration
$Z$ is a $(d,m,a)$-{\bf candidate} if the following two conditions hold:
$$h^0(X, \fasm O(d))\leq \mbox{deg} (Z)$$
and
$$\mbox{deg}\left( Tr_{G_a}(Z)\right)\leq h^0(G_a,\fasm
O_{G_a}(d)).$$
We consider a $(d,m,a)$-candidate $Z$ to be a candidate for the
property
$h^0(X,\fasm I_Z(d))=0$ and we say that $Z$ is {\bf winning} if this
property holds.\end{defn}
The bound for $\mbox{deg} (Z)$  in the definition of candidates
is for convenience: a vanishing statement for a more general
configuration will be treated by considering the candidate obtained
by adding the right number of simple points.

The following easy lemma says that for large $d$, candidates contain
sufficiently many free points.
\begin{lem}\label{free} Let $m$ and $a$ be positive integers. For any
$(d,m,a)$-candidate $Z$, we have
$$\mbox{deg}(\mbox{Free}(Z))\geq \nu
\frap{d^{n+1}}{(n+1)!}- O(d^n),$$
where $\nu$ is the degree of $X$.
\end{lem}
For presentation purposes we introduce the
\begin{defn}\label{dvm} Given a polarised pair $(V,\fasm O(1))$ and
$m>0$ we define $\mbox{\bf d}(V,m)$ to be the least degree
(a-priori possibly infinite, and a-posteriori finite by our theorem)
such that for $d\geq \mbox{\bf d}(V,m)$ any
$(d,m,0)$-candidate is winning.\end{defn}
%

%
\section{Derivatives}
%
In practice, when we apply lemma \ref{diff}, we think of the
condition (dime) as being satisfied. This is easily
justified by an induction hypothesis on the
dimension (i.e. precisely that $d(G_a,m)$ is finite).
Lemma \ref{diff} is then a justification for replacing $T$ by
$T_{\star}\up$.

This leads us to introduce a formal operator $\gotm D$ sending
one $(a,m)$-configuration to another which we call the derivative
(see \ref{D1}). Of course, we are especially interested in the case
where this operator takes $(d,m,a)$-candidates to
$(d-a,m,a)$-candidates. In the present section, we define the
derivative and show that it behaves well for large $d$.

Here is the idea behind the definition of a derivative. Given
a $(d,m,a)$-candidate $Z$, we wish to apply our lemma \ref{diff}
as follows. We specialize the $s$ biggest free points of
$Z$ onto the divisor $G_a$, with $s$ as large as possible. Still a few
conditions (say $r$) are missing in $G_a$, and we require that $r$
further free points be available in $Z$ so that we may apply
\ref{diff}. In that case, the derivative of $Z$  is the subscheme
$T\up_{\star}$ involved in the degue condition of \ref{diff}.
\begin{defn}\label{D1}
Let $Z$ be a $(d,m,a)$-candidate on $X$ with $t=t(Z)$ free points
$P_1^{m_1},\ldots ,P_t^{m_t}$, where the multiplicities
appear in non-decreasing order. Let $s=s(Z)\leq t$ be
the greatest integer such that
$$\mbox{deg}
(Tr_{G_a}(Z))+\mbox{$\binom{m_1+n-1}{n}+
\cdots +  \binom{m_s+n-1}{n}$}\leq h^0(G_a,\fasm
O_{G_a}(d)).$$
and set
$$
r=r(Z) = h^0(G_a,\fasm O_{G_a}(d)) -
\mbox{deg}
(Tr_{G_a}(Z))-\mbox{$\binom{m_1+n-1}{n}-
\cdots -  \binom{m_s+n-1}{n}$}.$$
We say that $Z$ is {\bf derivable} with respect to $G_a$ if
$$
r+s\leq t.
$$
If $Z$ is a derivable $(d,m,a)$-candidate, its {\bf derivative} with respect to
$G_a$, denoted $\gotm D^{(1)}(Z)$, is defined to be the $(a,m)$-configuration
$$\br{lll} \gotm D^{(1)}(Z)&=&P_{s+r+1}^{m_{s+r+1}}\; \cup\;
\cdots\;
\cup \; P_t^{m_t}\; \cup\;  \mbox{Const}(Z)\up\; \cup \\
&&\\
&&\hspace{3cm} Q_1^{m_1-1}\; \cup
\;\cdots\;  \cup \; Q_s^{m_s-1}\; \cup\\
&&\\
&&\hspace{4.5cm}
 D_{m_{s+1}}(Q_{s+1})\; \cup\;
\cdots\;  \cup\;  D_{m_{s+r}}(Q_{s+r})
\er$$
where $Q_1,\ldots ,Q_{s+r}$ are generic points of $G_a$ and the
notation is that of \ref{trace-residual} and \ref{simple residue}.
\end{defn}
Recall that $\alpha_0$ is an integer such that $\fasm O(a)$
is very ample for $a \ge \alpha_0$. What we need
to know about the derivative is the following :
\begin{lem}\label{der}  Let $a\geq \alpha_0$ and $m$ be positive
integers. Then there exists an integer $\mbox{\bf der}(a,m)$ such that
for any $d\geq \mbox{\bf der}(a,m)$ and any $(d,m,a)$-candidate $Z$ on
$X$:
\be
\item
$Z$ admits a derivative $\gotm D^{(1)}(Z)$;
\item
for any $N$, if $Z$ has either no free point of multiplicity $m$ or
at least $N$ free points of multiplicity less than $m$, then so does
$\gotm D^{(1)}(Z)$;
\item
The degree of the trace of $\gotm D^{(1)}(Z)$
satisfies the following estimate, where, as above,
$\nu = deg (\fasm O(1))$:\\
%
$\begin{array}{lll}
\mbox{deg Tr}_{G_a}\left(
\gotm D^{(1)} (Z)\right)
&\;=\;&
\left(\frap{(m-1)\, a\, \nu}{m+n-1}\right)\;
\frap{d^n}{n!}+O(d^{n-1})\\
&\;=\;&
h^0(G_a,\fasm O_{G_a}(d-a))-
\left(\frap{n\, a\, \nu}{m+n-1}\right)\;
\frap{d^n}{n!}+O(d^{n-1});\\
\end{array}$

\item
$\gotm D^{(1)}(Z)$ is
a $(d-a,m,a)$-candidate;
\item
if $d(G_a,m)$ is finite and  $\gotm D^{(1)}(Z)$ is winning, then so
is $Z$.
\ee
\end{lem}

\noindent{\bf Proof.}

In order to prove 1.,  it is enough to prove that the number of free points in
$Z$ is larger than $2h^0\bigl(G_a,\fasm O_{G_a}(d)\bigr)$. The latter
is bounded by $Cd^n$ for some constant $C$, so we may conclude by \ref
{free}.

As for 2., it is an immediate consequence of the definition of the
derivative.

For 3., let $r,s, t$ and $m_i$ be as in \ref{D1}. Then
$$\frap{m+n-1}{n}\left(\mbox{$\sum_{i=1}^{s}$}\mbox{$\binom{m_i
+n-2}{n-1}$}\right)\geq
\mbox{$\sum_{i=1}^{s}$}\mbox{$\binom{m_i
+n-1}{n}$}=h^0(G_a,\fasm O_{G_a}(d))-r$$
and
$$\br{ccl}\mbox{deg Tr}_{G_a}\left(\gotm D^{(1)}(Z)\right)&\leq&
\sum_{i=1}^{s}\mbox{$\binom{m_i +n-2}{n}$}+r\,
\mbox{$\binom{m+n-1}{n}$}
\vspace{1ex}\\
& = & \left( \sum_{i=1}^{s} \mbox{$\binom{m_i+n-1}{n}$}+r\right) -
\sum_{i=1}^{s}\mbox{$\binom{m_i+n-2}{n-1}$}\vspace{1ex}\\
&&- r+r
\,\mbox{$\binom{m+n-1}{n}$}
\vspace{1ex}\\
& \leq & h^0(G_a,\fasm O_{G_a}(d))
-\frap{m-1}{m+n-1}\biggl( h^0 (G_a,\fasm O_{G_a}(d)) -
r)\biggr)\vspace{1ex}\\
&&- r +r\, \mbox{$\binom{m+n-1}{n}$}\vspace{.5ex}\\
& \leq & \frap{n}{m-1}\left( h^0 (G_a,\fasm
O_{G_a}(d))\right)+\mbox{$\binom{m+n-1}{n}$}^2\vspace{1ex}\\
&= & \left(\frap{(m-1)\, a\,
\nu}{m+n-1}\right)\; \frap{d^n}{n!}+O(d^{n-1}).
\er$$
Finally, we have
$$\br{ccl}
h^0(G_a,\fasm O_{G_a}(d-a))-\mbox{deg
Tr}_{G_a}\left(\gotm D^{(1)}(Z)\right) & \geq & h^0(G_a,\fasm
O_{G_a}(d))-\frap{m-1}{m+n-1}h^0(G_a,\fasm O_{G_a}(d))
\vspace{1ex}\\&&-\mbox{$\binom{m+n-1}{n}^2$}
\vspace{1ex}\\
&= &
\left(\frap{n\, a\, \nu}{m+n-1}\right)\;
\frap{d^n}{n!}+O(d^{n-1}).
\er$$

For 4., we first note that, when $\gotm D^{(1)}(Z)$ is
defined and  $H^1(X,\fasm O(d-a))=0$, one has
$$h^0(X,\fasm O(d-a))\leq \mbox{deg}\left(\gotm D^{(1)}(Z)\right) =
\mbox{deg}(Z)-h^0(G_a,\fasm O_{G_a}(d)).$$
This means that,
for sufficiently large $d$,
the $(a,m)$-configuration $\gotm D^{(1)}(Z)$ is a
$(d-a,m,a)$-candidate, since by 3., its trace on $G_a$ has degree at most
$h^0(G_a,\fasm O_{G_a}(d-a))$.

For 5., using the notation of \ref{D1}, we apply \ref{diff}, with $Z_0$ the
closed subscheme
$$\mbox{Const}(Z)\; \cup \; Q_1^{m_1}\; \cup \cdots \cup \;
Q_s^{m_s}\;\cup \; P_{s+r+1}^{m_{s+r+1}}\; \cup \cdots \cup\;
P_t^{m_t}.$$
Let $W=Q_1^{m_{s+1}}\cup \cdots
\cup Q_r^{m_{s+r}}$.
The dime of \ref{diff} holds for $d\geq d(G_a,m)$,
while the degue of \ref{diff} is just the hypothesis that
$\gotm D^{(1)}(Z)$ is winning, so the lemma follows from \ref{diff}.
$\hfill \square$
%

%
\section{Concentrated derivatives} \label {CD}
%

If theorem \ref{thmm} were true for low degrees, then
repeated applications of lemma \ref{diff}, hence of the derivative,
would suffice to prove the theorem by induction on the degree. Instead
one must modify the process and try to reduce the multiplicities of
the free points, thus ending the proof by induction on the
multiplicity. This is done using a specialisation of the second
derivative (see \ref{D(2)}): bearing  in mind the semi-continuity of
the cohomology, one easily sees that the (degue) of \ref{diff} holds
if it holds for some specialisation of $T\up_{\star}$. A complication
arises with the degree of the base divisor $G_a$ which must
be lowered during the induction on $d$ before an induction hypothesis
on $m$ allows one to finish the proof. We get around this problem using
a specialisation of the first derivative which we call a concentrated
derivative. In this section we introduce this concentrated
derivative and prove  results analogous to those for derivatives.

\begin{defn}\label{defcd} Let $d\, ,m\, , a$ be positive integers with
$a>1$, and let
$Z$ be a derivable $(d,m,a)$-candidate. We define the {\bf
concentrated derivative} of $Z$ with respect to $G_a$, denoted
$\gotm D^{(1)}_c(Z)$, to be the
$(a-1,m)$-configuration obtained from $\gotm D^{(1)}(Z)$
by degenerating $G_a$ to the
generic union $G_1+ G_{a-1}$ and specialising all $G_a$-residues of
$\gotm D^{(1)}(Z)$ to have generic support in $G_{a-1}$.
\end{defn}
What we need to know about the concentrated derivative is concentrated in
the following:
\begin{lem}\label{derc}
Given $m>0$
there exists an integer $A(m)$ such that for all $a\geq
A(m)$ there exists an integer $\mbox{\bf derc}(a,m)$ such that for any
$d\geq
\mbox{\bf derc}(a,m)$ and any $(d,m,a)$-candidate $Z$:
\be
\item
$Z$ admits a
concentrated derivative $\gotm D^{(1)}_c(Z)$ which is
a $(d-a,m,a-1)$-candidate;
\item
for any $N$, if $Z$ has either no free point of multiplicity $m$ or
at least $N$ free points of multiplicity less than $m$, then so does
$\gotm D^{(1)}_c(Z)$;
\item
if $d(G_a,m)$ is finite and  $\gotm D^{(1)}_c(Z)$ is winning, then so
is
$Z$.
\ee
\end{lem}

\noindent{\bf Proof.}

For 1., let $A(m)$ be an integer
$a$ satisfying
$$\frap{m-1}{m+n-1}A(m) < A(m)-1.$$
Then, for $a\geq A(m)$ and $d$ sufficiently large,
for any $(d,m,a)$-candidate $Z$,  we have, by (\ref{der}.3),
$$\mbox{deg}\, \mbox{Tr}_{G_{a-1}}\left( \gotm D^{(1)}_c(Z)\right)
= \mbox{deg}\,
\mbox{Tr}_{G_a}\left( \gotm D^{(1)}(Z) \right)\leq h^0(G_{a-1},
\fasm O_{G_{a-1}}(d-a))$$
so that $\gotm D^{(1)}_c(Z)$ is a $(d-a,m,a-1)$-candidate.

As for 2., it follows from the similar statement for the
derivative, since the derivative and the concentrated
derivative have the same free points.

For 3., if $\gotm D_c^{(1)}(Z)$
is winning then so is $\gotm D^{(1)}(Z)$, since the former is a
specialisation of the latter. We conclude that $Z$ is winning for
$d\geq \mbox{\bf der}(a,m)$ by (\ref{der}.5).$\hfill \square$

\section {Special second derivative}

In this section, we explain the construction which generates free
points. This corresponds to a modified second derivative, which we
denote by $\gotm D^{(2)}[\alpha]$, where $\alpha$ is an integer.
\begin{defn}\label{D(2)}Let $m,a>0 $, and let $Z$ be a  twice
derivable $(d,m,a)$-candidate. Let $r^{(2)}(Z)$ be the number of
residues of $\gotm D^{(2)}(Z)$ which are points, necessarily of
multiplicity at most $m-1$. For $0<\alpha <a$, we set
$$r^{(2)}[\alpha](Z)=\mbox{min}\left(h^0(X,\fasm
O(\alpha))-1,r^{(2)}(Z)\right)$$
and define $\gotm
D^{(2)}[\alpha](Z)$ to be the
specialisation of the second derivative $\gotm
D^{(2)}(Z)$ obtained by degenerating $G_a$ and its residues to the
generic union
$G_{\alpha}+G_{a-\alpha}$ with
$r^{(2)}[\alpha](Z)$
of the residues which are points specialised to
have generic support in $G_{\alpha}$, and all other residues
specialised to $G_{a-\alpha}$.\end{defn}
Here is what we need to know about this construction.

\begin{lem}\label{sder} Given $m,N$,
there exist integers $\alpha$ and
$a_0>\alpha $ such that for all $a\geq a_0$ there exists
$d'_0=d'_0(m,N,a)$ such that for $d\geq d'_0$ and any
$(d,m,a)$-candidate $Z$:
\be
\item
$Z$
is twice derivable and $\gotm D^{(2)}[\alpha](Z)$ is a $(d-2a , m
,a-\alpha)$-candidate having  either no free point of multiplicity $m$
or at least $N$ free points of multiplicity at most $m-1$;
\item
if  $d(G_a,m)$ and $d(G_{a-\alpha},m)$ are finite and
$\gotm D^{(2)}[\alpha](Z)$ is winning, then so
is
$Z$.
\ee
\end{lem}

\noindent{\bf Proof.}  Let $\alpha=\alpha(N)$ be an
integer satisfying $h^0(X,\fasm O(\alpha)) > N.$
For 1., let $a_0> \alpha  $ be  such that for
$a-\alpha > \frac{m-1}{m+n-1}a$ for $a\geq a_0$. Then for $a \geq a_0$,
and $d\vg 0$, we have
$$\mbox{deg}\, \mbox{Tr}_{G_a}\left(\gotm D^{(2)}(Z)\right)
< h^0(G_{a-\alpha}, \fasm O_{G_{a-\alpha}}(d-2a))$$
because by (\ref{der}.3)
$$\mbox{deg}\, \mbox{Tr}_{G_{a}}\left(\gotm
D^{(2)}(Z)\right)\leq
\frap{m-1}{m+n-1}a\nu
\frap{d^n}{n!}+ O(d^{n-1})$$
while
$$h^0(G_{a-\alpha} , \fasm O_{G_{a-\alpha}}(d-2a))\geq
(a-\alpha)\nu\frap{d^n}{n!}+ O(d^{n-1}).$$
This implies that $\gotm
D^{(2)}[\alpha](Z)$ is a $(d-2a , m ,a-\alpha)$-candidate.

Now for $d\vg 0$, by (\ref{der}.3), we have
$$h^0(G_{a},\fasm O_{G_a}(d-a))-\mbox{deg}\, \mbox{Tr}_{G_a}\,
(\gotm D^{(1)}(Z))\geq N\mbox{$\binom{m+n-1}{n}$}$$
so that, in the notation of \ref{D1}, $s(\gotm D^{(1)}(Z))\geq N$.
If $\gotm D^{(1)}(Z)$ has at least $N$ free points of multiplicity
$m$, then $\gotm
D^{(2)}[\alpha](Z)$ has $N$ free points of multiplicity $m-1$:
indeed, in that case,
\mbox{$r^{(2)}[\alpha](Z)$} is larger than $N$, and the
\mbox{$r^{(2)}[\alpha](Z)$}
points specialised to
\mbox{$G_{\alpha}$} are without constraint, since any set of
$N$ points lie on an effective divisor in the
linear system \mbox{$\bigl|H^0(X,\fasm O(\alpha))\bigr|$}.
Otherwise,
$\gotm D^{(2)}[\alpha](Z)$ has no more free points of  multiplicity $m$.

For 2., we observe that the second derivative $\gotm D^{(2)}(Z)$
is also a winning candidate, and conclude by applying twice (\ref
{der}.5). $\hfill \square$
%

\section{Proof of the theorem}
%
\subsection{A proposition implying the theorem}
The following proposition (which we prove below \ref{prp}) sums up the
efforts of the previous sections and, as we willl now show, easily
implies our theorem.

\begin{prop}\label {mainp} Let $X$ be a projective, geometrically
reduced and irreducible variety of dimension $n+1$ over a field $k$
of arbitrary characteristic $p$. If $p > 0$,
suppose further that $X$ is smooth in
codimension one if $p>0$. Let $\fasm O(1)$ be an invertible ample
bundle on $X$. Given $m>0$, there exists $a_0(m)$ such that for any
$a\geq a_0(m)$ there exists $d_0(a,m)$ such that for all
$d\geq d_0(a,m)$ any $(d,m,a)$-candidate is winning.\end{prop}
\noindent{\bf Proof of \ref {thmm}.}

We first handle the case where $\fasm M = \fasm O$.

We take $d_0(m)=d_0(a_0(m),m)$ and consider some
$d\geq d_0$ and some generic union $Z$ of (fat) points of
multiplicity $\leq m$. If the degree of $Z$ is smaller than
$h^0(\fasm O(d))$, we reduce to the case with equality by adding
generic simple points. Since the trace on $G_{a_0}$ is empty, we may
consider $Z$ as a $(d,m,a_0)$-candidate (\ref{candidates}), and
conclude by \ref {mainp}.

As announced, we only gloss over the proof in the case where $\fasm M$
is  arbitrary.

Firstly, replacing $\fasm M$ by $\fasm M \otimes \fasm L ^b$,
we may suppose that $\fasm M$ is effective. Next, we can suppose as above
$deg \, Z \geq h^0(\fasm M \otimes \fasm L^d)$ and we have
to prove that
$H^0(X, \fasm I_Z \otimes \fasm M \otimes \fasm L^d)=0$.
The idea of the proof is then to choose a suitable $a$ and
to apply \ref {diff} as in \ref{der} using the generic divisor $G^{\star}_a$ in
$\vert H^0(X,\fasm M \otimes \fasm L ^a) \vert$. By induction on the
dimension we
can suppose that the dime condition holds. To prove the degue condition, we
degenerate
$G^{\star}_a$ to $M \cap G_a$, where $M$ is in $\vert H^0(X,\fasm M )\vert$ and
$G_a$ is the generic divisor in $\vert H^0(X, \fasm L^a) ,\vert$;
specializing all
the residues onto $G_a$. In this way, we get a $(a,m)$-configuration
$Z\up_c$ and,
if this is a $(d-a, m, a)$-candidate, we can end with the particular case
($\fasm M= \fasm O$) since such a candidate is winning
for sufficiently large $d$.
To see that $Z_c^{\prime}$ is a $(d-a,m,a)$-candidate
for suitable large $a$ and $d$, one persues
an argument analogous to \ref{der}.1 and one shows that
for any
$a$, and all $d$ sufficiently large with respect to $a$, $Z$ has enough free
points to make \ref {diff} applicable. Then, as in \ref{derc}.1, one shows that
for sufficiently large $a$ and all $d$ sufficiently large with respect to
$a$, the
$(a,m)$-configuration $Z\up_c$ is a $(d-a, m, a)$-candidate.$\hfill \square$
\begin {rem}
A further generalisation would be to take a fixed closed (zero-dimensional)
subscheme $V_0$ and its union with points of
multiplicity $\leq m$. The union of $V_0$ with sufficiently many
generic simple points has maximal rank in all degrees giving the
initial case for an induction on the multiplicity, while the proof of
the dimension one case is virtually unchanged.
\end{rem}
\subsection{\label{prp}Proof of the proposition}
To prove the proposition, we argue by induction on the dimension
$n+1$. Note that in all characteristics, the generic effective divisor
in a very ample linear system on a variety $X$ of dimension $>1$ is a
variety which is smooth outside the singular locus of $X$ (see [L]
VII 13). Thanks to the initial cases \ref{ci} and \ref{cpi} below, we
may suppose that the proposition has been proven for multiplicity $m$
in dimension $n$ and for multiplicity $m-1$ in dimension $n+1$. This
implies that
$\mbox{\bf d}(G_a,m)$ is finite for all $a\geq 1$ and that there
exists $a_0(m-1)$ such that for $a\geq a_0(m-1)$ there exists
$d_0(a,m-1)$ such that for $d\geq d_0(a,m-1)$ any
$(d,m-1,a)$-candidate is winning.  We proceed in three steps.

\noindent{\bf First step.} With the notation of \ref{der} and
\ref{derc}, we define
$$b_0=\mbox{max}\bigl(A(m), a_0(m-1)\bigr),$$
$$\Delta=\mbox{max}\bigl(\mbox{\bf
der}(b_0,m),d_0(b_0,m-1)+b_0 \bigr)$$
and
$$N=h^0(X,\fasm O(\Delta+b_0 -1)) + \binom{n+m-1}{n},$$
and prove by induction that, for any $d \geq \Delta$, any
$(d,m,b_0)$-candidate with either no free point of multiplicity $m$ or
at least $N$ free points of multiplicity less than $m$ is winning.

We start with the case $\Delta \leq d < \Delta + b_0$,
and consider a $(d,m,b_0)$-candidate $Z$ with
either no free point of multiplicity $m$ or
at least $N$ free points of multiplicity less than $m$.
If $\mbox{deg}(Z)\geq h^0(X,\fasm O(d)) + \binom{n+m-1}{n}$,
and $Z$ has a free point of multiplicity $m$,
we may replace $Z$ by the subscheme obtained by diminishing
by one the multiplicity of this free point, which still has
at least $N$ free points of multiplicity less than $m$.
In other words, we may suppose either
that $Z$ has no free point of multiplicity $m$,
or that $\mbox{deg}(Z) < h^0(X,\fasm O(d)) + \binom{n+m-1}{n}$
holds. In the latter case,
there is no room for $N$ free points of multiplicity less than $m$.
Summing up, we can suppose that $Z$ has no free point of multiplicity $m$.
Thanks to $d \geq \mbox {\bf der}(b_0,m)$ and \ref {der}, $Z$ has a
first derivative $\gotm D^{(1) }(Z)$  which is a
$(d-b_0,m-1,b_0)$-candidate. Thanks to $d-b_0 \geq d_0(b_0,m-1)$,
this candidate is winning. Thanks to $d \geq \mbox {\bf der}(b_0,m)$
and \ref {der} again,
$Z$ is winning too.

For $d\geq
\Delta + b_0$ let $Z$ be a $(d,m,b_0)$-candidate having
either no free point of multiplicity $m$, or
at least
$N$ free points of multiplicity at most $m-1$.
Thanks to $d \geq \mbox {\bf der}(b_0,m)$
and \ref {der}, $Z$ has a first derivative $\gotm D^{(1) }(Z)$
which is a $(d-b_0,m,b_0)$-candidate having
either no free point of multiplicity $m$ or
at least
$N$ free points of multiplicity at most $m-1$. Thanks to the inductive
assumption, $\gotm D^{(1) }(Z)$ is winning. Again thanks to
$d \geq \mbox{\bf der}(b_0,m)$ and \ref {der}, $Z$ is winning too.

\noindent{\bf Second step.} Here we prove that for any  $b\geq b_0$ there
exists $\delta =\delta(b,m)$ such that for $d\geq \delta $ any  $(d,m,
b)$-candidate having either no free point of multiplicity $m$, or
at least $N$ free points of multiplicity at most $m-1$ is winning.

The proof is by induction on $b$. The initial case $b=b_0$ is the
previous step. For the induction step,
we take
$\delta (b) = max ( \mbox {\bf derc}(b,m), \delta (b-1)+b)$. The
statement then follows by \ref {derc}, which applies because  $b_0 \geq
A(m)$.

\noindent{\bf Final step.} Here we set $a_0=a_0(m)=max (b_0+\alpha(N),
a_0 (m,N))$ where $\alpha=\alpha(N)$ and $a_0 (m,N))$ are defined in \ref
{sder}, and, for $a \geq a_0$,  $d_0=d_0(a,m)= max (d'_0 (m,N,a),
\delta(a-\alpha,m) +2a)$, and we prove the full statement, namely that,
for $d\geq d_0(a,m)$, any $(d,m,a)$-candidate $Z$ is winning.

Indeed, by \ref {sder} applied to $n,N$,  $Z$
is twice derivable and $\gotm D^{(2)}[\alpha](Z)$ is a $(d-2a , m
,a-\alpha)$-candidate having  either no free point of multiplicity $m$
or at least $N$ free points of multiplicity at most $m-1$.
Since $d-2a \geq \delta(a-\alpha,m)$, this candidate is winning
by the second step. This implies that $Z$  itself is winning by
\ref {sder}. $\hfill \square$
%

\subsection{The proposition in dimension one}
%
The initial case $n=0$ can be deduced from the following
general results for curves. We first treat the characteristic zero
case with the
\begin{prop} \label{ci}
Let $C$ be a geometrically irreducible
quasi-projective curve over a field $k$ of characteristic zero. Let
$V\subset H^0(C,\fasm L)$ be a linear subspace of finite dimension $v$
of global sections of the invertible sheaf $\fasm L$ on $C$. Let
$x_1,\ldots ,x_r$ be the generic set of $r$ closed points of $C$
defined over the function field $K$ of
$C\times \cdots \times C$ ($r$ factors), and let
$m_1,\ldots ,m_r$ be positive integers. Let $D$ be the divisor $m_1x_1
+ \cdots + m_rx_r$ on $C_K=C\times_k K$. Then the canonical map
$$V\lra H^0(C_K,\fasm O_D\otimes \fasm L)$$
has maximal rank.\end{prop}

\noindent{\bf Proof.} If $v\neq m=\sum_i m_i$ , one can either
diminish the multiplicities or add (generic) free points  and
suppose that $v=m$. Since the property is open, we can specialise to
the case of a single point $x$ and the divisor $D=mx$. In
this case the proposition is equivalent to showing that the
determinant of the canonical map
\begin{equation}\label{pp}V\otimes \fasm O_C \lra \mbox{P}^v(\fasm
L)\end{equation} is not identically zero, where $\mbox{P}^v(\fasm L)$
is the sheaf of $v^{th}$ order principal parts of $\fasm L$. For this
we can suppose that the base field is algebraically closed and, since
this map commutes with localisation and the completion at a closed
point of $C$, it is sufficient to show that the canonical map
$$V\otimes k[[t]]\lra k[[t,x]]/((x-t)^{v})$$
$$ f \mapsto
f(t)+f^{\prime}(t)(x-t)+ f^{\prime\prime}(t)\frac{(x-t)^2}{2!} +\cdots
+f^{(v-1)}(t)\frac{(x-t)^{v-1}}{(v-1)!}$$
has maximal rank. Choosing a basis $f_1,\ldots ,f_v$ for $V$, the
determinant of this map is just the Wronskian
$$W(f_1,\cdots ,f_v)=\mbox{det}\left[\frac{\partial^i
f_j}{\partial t^{i}}\right]$$
which, as is well known, has maximal rank for $f_1,\ldots ,f_v$
linearly independent.$\hfill \square$

We now give the initial case for smooth curves in arbitrary
characteristic.
\begin{prop}\label{cpi} Let $C$ be a smooth, geometrically connected,
projective
curve of genus $g$ over an arbitrary field. Let $\fasm M$, $\fasm
L$ be line bundles on $C$ with $\fasm L$ ample, let $m>0$ be an
integer and let $d_0(m)$ be the least integer $d$ such that
$h^0(C,\fasm M\otimes \fasm L^d) > m(m-1)(g-1)/2$ and $\fasm M\otimes
\fasm L^d$ is non-special. Let
$x_1,\ldots ,x_r$ be generic points on
$C$ and let $Z$ be the divisor $m_1x_1+\cdots +m_rx_r$ where
$0<m_i\leq m$ for $i=1,\ldots ,r$. Then the canonical map
$$H^0(C,\fasm M\otimes \fasm L^d)\lra H^0(C,\fasm O_Z\otimes
\fasm M\otimes \fasm L^d)$$
has maximal rank for $d\geq d_0(m)$.
\end{prop}

\noindent{\bf Proof.} Adding points if necessary, we can suppose
$\mbox{deg}(Z)\geq h^0(C,\fasm M\otimes \fasm L^d)$. By hypothesis we
then have $\mbox{deg}\, Z > m(m-1)(g-1)/2$ so that some set of $g$
points amongst the $x_i$ have the same multiplicity $m_0$. Renumbering,
we can write $m_1x_1+\cdots + m_rx_r =m_0(y_1+\cdots +y_g) + D=Z+D$, where
$D$ has support away from the $y_i$. Since the natural map $C^g\lra
\mbox{Pic}^g(C)$ and the power map
$\mbox{Pic}^g(C)\lra \mbox{Pic}^{m_0g}(C)$ are surjective, it follows
that for $y_1,\ldots ,y_g$ generic, the sheaf $\fasm O(Z)$ and hence $\fasm
L^d\otimes
\fasm M\otimes \fasm O(-Z)$
is the generic sheaf in   its component of the Picard
scheme so that either
$h^0(C,\fasm L^d\otimes \fasm M\otimes \fasm O(-Z))=0$ or $h^1(C,\fasm
L^d\otimes \fasm M\otimes \fasm O(-Z))=0$. \hfill$\hfill \square$

This completes the proof of the cases in dimension 1. We end with the
following example showing that the `smooth in codimension one'
hypothesis cannot be dropped in characteristic $p>0$.

\begin{rem}\label{cp} Let $p$ be an odd prime and
$C$ the plane curve defined by
the equation $y^2-x^p=0$ over an
algebraically closed field of characteristic $p$. The tangent line at
$z=(t^2,t^p)$, $t\neq 0$, is given by $y=t^p$ and has a contact of order
$p$ with $C$ at $z$. It follows that for any choice $z_1,\ldots ,z_d$
of points on the smooth locus of $C$, the divisor $Z=pz_1 +\cdots
+pz_d$ is an effective divisor associated to $\fasm O_C(d)$, whereas
$h^0(C,\fasm O_C(d))= dp+1-(p-1)(p-2)/2\leq dp$ for $d\geq p-2$.
\end{rem}
%
\section{\label{formallemma}The formal lemma}

In this section, we prove the formal part of our
differential lemma, the rest of the proof being in the next section.
We would like to point out that the original motivation and proof of the
following results owed much to the work [LL1,2].
\subsection{Preliminaries}
Consider the algebra of formal functions $k[[\bm x,y]]$, where $\bm
x=(x_1,\ldots ,x_{n-1})$, which we furnish with an ideal $I$ of the form
$$I=I_0\oplus I_1y\oplus \cdots \oplus I_{m-1}y^{m-1}\oplus (y^m)$$
where, for $\alpha=0,..., m-1$,
$I_{\alpha}\subset k[[\bm x]]$ is an ideal. We call such ideals
{\em vertically graded ideals}. Note that
\begin{equation}\label{inclusions}I_0 \subset I_1 \subset \cdots
\subset I_{m-1}\end{equation}
An ideal
$$I_t=I_0[[t]]\oplus I_1[[t]](y-t^r)\oplus \cdots\oplus
I_{m-1}[[t]](y-t^r)^{m-1}\oplus ((y-t^r)^m)$$
in the algebra $k[[t,\bm x,y]]$ is called a {\em standard
deformation } of the vertically graded ideal $I$. For $i\geq m$ we
let $I_i=k[[\bm x]]$.

Given  a function $F_0+F_1t+\cdots $ in $I_t$, the functions $F_i(\bm x
,y)$ must satisfy certain residual conditions. If $r=1$ and $I=(\bm
x,y)^m$, the residual condition is just that $F_i(\bm x,y)$ must
vanish to the order $m-i$, and can be compared with [Xu]. This is the
sense of the following statement.
\begin{prop}\label{gettingstarted} Let $F =\sum_{\alpha\geq
0}F_{\alpha}(\bm x,y)t^{\alpha}
=\sum_{\alpha, \beta\geq 0} F_{\alpha,\beta}(\bm x)t^{\alpha}y^{\beta}$
be a function in $I_t$.
Then
$$\br{lll}F_{\alpha,\beta}(\bm x)
&\in &I_{\beta+[\![
\frac{\alpha}{r}]\!]}\er.$$
If $y$ divides $F_{\alpha}$ for $\alpha = 0,r,2r,\ldots ,pr$
then $F_0(\bm x,y)$ is in the ideal
$$I_0y\oplus I_1y^2\oplus \cdots \oplus I_{p-1} y^{p}\oplus
I_{p+1}y^{p+1}
\oplus\cdots \oplus I_{m-1}y^{m-1}\oplus ((y^m)).$$
\end{prop}

\noindent{\bf Proof.} Write $F$ in the following form
$$F = a_0(\bm x,t) + a_1(\bm x,t)(y-t^r)+\cdots +a_{m-1}(\bm
x,t)(y-t^r)^{m-1} + a_m(\bm x,t)(y-t^r)^m + \cdots $$
with
$$a_i(\bm x,t)=\sum_{j\geq 0}a_{ij}(\bm x)t^j$$
hence $a_{ij}(\bm x)\in I_i$. Developping out we find
$$\br{lll}F_{\alpha ,\beta}&=&\sum_{\nu =0}^{[\![
\frac{\alpha}{r}]\!]} \; (-1)^{\nu}\, \binom{\beta+\nu}{\beta}
a_{\beta+\nu,\, \alpha-\nu r}(\bm x)\\
&\in  &I_{\beta+[\![
\frac{\alpha}{r}]\!]}\er$$
where $[\![ z]\!]$ is the greatest integer part of $z$. This
proves the first part.

Now suppose that $y$ divides $F_{\alpha}$ for $\alpha= \lambda r$ and
 $\lambda=0,1,\ldots ,p$. Then we have
$$0=F_{\lambda r,\, 0} = a_{0,\lambda r}-a_{1,(\lambda - 1)r}+\cdots
+(-1)^{\lambda -1}a_{\lambda -1,r}+ (-1)^{\lambda}a_{\lambda,0}$$
so that $a_{0,0}=0$ and $a_{\lambda,0}\in I_{\lambda -1}$
for $\lambda =1,\ldots ,p$ as one sees using $a_{\mu , \nu}\in
I_{\mu}$ and (\ref{inclusions}). This gives the last part of the
proposition.  $\hfill \square$
\subsection{The formal lemma}
Throughout this subsection we will use the following notation.

For $i=1,\ldots ,\ell$, let $B^{(i)}=k[[\bm x_i,y_i]]$ be
an algebra of formal functions in $n$ variables where
$\bm x_i=(x_{i,1},\ldots ,x_{i,n-1})$ and let
$$I^{(i)}=I_0^{(i)}\oplus I_1^{(i)}y_i \oplus \cdots\oplus
I_{m_i-1}^{(i)}y_i^{m_i-1}\oplus (y_i^{m_i})$$
be a vertically graded ideal in $B^{(i)}$. Let
$$I=I^{(1)}\times \cdots \times I^{(\ell)}\subset B^{(1)}\times \cdots
\times B^{(\ell)}=B.$$
Let $k[[\bm t]]=k[[t_1,\ldots ,t_{\ell}]]$ and let $I_{\bm t}$ in
$B[[\bm t]]$ be the product of the ideals
$$I^{(i)}_{\bm t}=I_0^{(i)}[[\bm t]]\oplus I_1^{(i)}[[\bm t]](y_i-t_i)
\oplus
\cdots \oplus I_{m_i-1}^{(i)}[[\bm t]](y_i-t_i)^{m_i-1}\oplus
((y_i-t_i)^{m_i}).$$

Let $y=(y_1,\ldots ,y_{\ell})$ and for any linear subspace $V\subset
B$, let $V_{\mbox{res}(y)}=\{ v\in B\, |\, vy\in V\,\}$. Since $y$ is
not a zero-divisor, we get a residual exact sequence
\begin{equation}\label{res.of.v}0\lra
V_{\mbox{res}(y)}\stackrel{y}{\lra}V\lra V/V\cap (y)\lra
0\end{equation}

\begin{prop}\label{fl} Let $V\subset B$ be a $k$-linear subspace.
Suppose that for $i=1,\ldots ,\ell$ there exist
nonnegative integers $p_i$
such that the following two conditions are satisfied
\be
\item the canonical map
$$ V/V\cap (y)\lra k[[\bm x_1]]/I^{(1)}_{p_1}\times \cdots \times
k[[\bm x_{\ell}]]/I^{(\ell)}_{p_{\ell}}$$
is injective
\item The canonical map
$$V_{\mbox{res}(y)}\lra B/J$$
is injective where $J=J^{(1)}\times \cdots \times J^{(\ell)}$ and
$$J^{(i)}=I^{(i)}_{0}\oplus I^{(i)}_{1}y_i\oplus \cdots
\oplus I^{(i)}_{p_i-1}y_i^{p_i-1}\oplus I^{(i)}_{p_i+1}y_i^{p_i}\oplus
\cdots \oplus I^{(i)}_{m_i-1}y_i^{m_i-2}\oplus (y_i^{m_i-1})$$
\ee
Then the canonical map
$$\varphi_{\bm t} : V\otimes k[[\bm t]]\lra B_{\bm t}/I_{\bm t}$$
is (generically) injective.
\end{prop}

\noindent{\bf Proof.} We first reduce to the
case where the
$p_i$ are positive.

Let us suppose for simplicity that $p_1,\ldots,p_s$ are positive
and $p_{s+1},\ldots ,p_{\ell}$ are all zero. We denote by $V_0$
the subspace of $V$ formed by the elements
vanishing in each of the $k[[\bm x_i,y_i]]/I^{(i)}$
for $i=s+1 , \ldots ,\ell$. Conditions 1. and 2.
of the proposition imply the corresponding conditions
for $V_0$ when only the first $s$ factors on the right
hand side are present.

If we write $\bm t'$ for $(t_1,\ldots,t_s)$, the conclusion of the
proposition in the case where all $p_i$ are positive then implies that
$V_0\otimes k[\bm t']]$ injects into $B_{\bm t'}/I_{\bm t'}$. Since
$\varphi_{\bm t}$ is a map of free $k[[\bm t]]$-modules,
it is enough to prove that its restriction $\varphi_{\bm t'}$ over
$Spec \, k[[\bm t']]$ is injective. We write
$$\varphi_{\bm t'}= (\varphi_{\bm t'}\up, \varphi_{\bm t'}\upp):
V\otimes k[\bm  t']]  \rightarrow B_{\bm t'}/I_{\bm t'} \times R$$
with
$$R= (B^{(s+1)}/I^{(s+1)} \times \cdots \times B^{(\ell )}/I^{(\ell
)}) \otimes k[[\bm t']]$$
The kernel of $\varphi_{\bm t'}\upp$ is $V_0\otimes
k[\bm t']]$ and the restriction of $\varphi_{\bm t'}\up$ to
this kernel is injective, thus so is $\varphi_{\bm t'}$.

Henceforth we suppose that the $p_i$ are positive and we let
$$h=lcm(p_1,\ldots ,p_{\ell})=r_ip_i$$
be the least common multiple of the $p_i$ and consider the
one-parameter  deformation obtained by setting $t_i=t^{r_i}$. Since the
rank of $\varphi_{\bm t}$ is semi-continuous, we need only show  that
the canonical map
$$\varphi_t : V\otimes k[[t]]\lra B[[t]]/I_t$$
obtained by the formal base change $k[[t_1,\ldots ,t_{\ell}]]\lra
k[[t]]; \; t_i\mapsto t^{r_i}\, $; is injective.

Let
$$F_t =(F^{(1)}_t,\ldots ,F^{(\ell)}_t)\in \ker\varphi_t\; =\; V_t\cap
I_t,$$
where $V_t$
is the image of $V\otimes k[[\bm t]]$ and $I_t$ is the image of $I_{\bm
t}$ in $B[[t]]$.

In case $F_0=0$, we may replace $F_t$ by $F_t/t$
since $B[[t]]/I_t$ is a torsion free $k[[t]]$-module, Thus we only have to
prove
$F_0=0$.

Since $F^{(i)}_t=\sum_{\alpha\geq 0}\, F^{(i)}_{\alpha}(\bm
x_i,y_i)t^{\alpha}\in I^{(i)}_t$, where $I^{(i)}_t$ is the image of
$I^{(i)}_{\bm t}$ in $B^{(i)}[[t]]$, the
first part of proposition
\ref{gettingstarted} implies
$$(F^{(1)}_{\alpha}(\bm x_1,0),\ldots ,F^{(\ell)}_{\alpha}(\bm
x_{\ell},0))\in I^{(1)}_{p_1}\times \cdots \times I^{(\ell)}_{p_{\ell}}$$
for $\alpha = 0,1,\ldots ,h$. Applying hypothesis 1. of the
proposition, we conclude that $y$ divides $(F^{(1)}_{\alpha}(\bm
x_1,y_1),\ldots ,F^{(\ell)}_{\alpha}(\bm x_{\ell},y_{\ell}))$ for
$\alpha = 0,1,\ldots ,h$.

Now applying the second part of proposition \ref{gettingstarted} we
obtain
$$F^{(i)}_0(\bm x_i,y_i) =y_i G_0^{(i)}(\bm x_i,y_i)$$
with $G_0^{(i)}(\bm x_i,y_i)\in  J^{(i)}$. Letting $G_0=(G_0^{(1)},\ldots
,G_0^{(r)})$ we see that $G_0$ is in $V_{\mbox{res}(y)}\cap J$, but the
second hypothesis of the proposition simply says that
$V_{\mbox{res}(y)}\cap J=0$, giving
> $F_0=0$ as required. $\hfill \square$
\section{The differential lemma}
Throughout this section,
 $X$ denotes an irreducible algebraic variety,
 $H$ a reduced irreducible positive Cartier divisor on $X$,
 $X^0$ a dense open nonsingular subscheme of $X$ such that
 $H^0:=H \cap X^0$ is the nonsingular locus of $H$.
 Finally $\fasm I_H$ denotes the ideal sheaf of $H$.

For $M$ another $k$-scheme,
we denote by $Hom(M, X)$ (resp. $Hom(M, X^0)$) the set of morphisms
from $M$ to $X$ (resp. $X^0$) as
well as, in case $M$ is projective, the corresponding Hilbert scheme.
If $M$ is algebraic, zero-dimensional and connected,
it is easy to check that the natural morphism from $Hom(M, X^0)$ to $X^0$
is smooth with smooth irreducible fibers. Thus
 $Hom(M, X^0)$ is also irreducible.
Its generic point represents an embedding whose image in $X$
we denote by $M_X$.

Now let $M$ be a subscheme of $\mbox{Spec}\, k[[x_1, \dots, x_n]]$.
We denote by $Hom(M, X, H)$ the set (or Hilbert scheme)
 of morphisms $f$ from $M$ to $X^0$
such that the ideal $f^*({\cal I}_H)$ is contained in $(x_n)$. We call these
morphisms $H$-morphisms from $M$ to $X$, and if a $H$-morphism is an
embedding, we say that it is a $H$-embedding.
 If $M$ is algebraic, thus zero-dimensional and connected,
it is easy to check that the natural
(restriction) morphism from $Hom(M, X, H)$ to $H^0$
is smooth. Furthermore, its fiber $Hom(M, X, H, z)$ over a
point $z$ of $H^0$ is a vector space, thus
smooth and irreducible. As a consequence, $Hom(M, X, H)$
is again irreducible and smooth.
Its generic point is a $H$-embedding whose image in $X$
we denote by $M_{X,H}$.

We say that the subscheme $M$ of $\mbox{Spec}\, k[[x_1, \dots, x_n]]$ is a
model
of dimension $n$ if its ideal is a vertically
graded  ideal as in \S 8:
$$I=I_0\oplus I_1x_n\oplus \cdots \oplus I_mx_n^m\oplus \dots$$
where
$I_m$ is a non-decreasing sequence of ideals
in $k[[x_1, \dots, x_{n-1}]]$ with
$I_m = k[[x_1, \dots, x_{n-1}]]$ for large $m$.

For $M$ a model of dimension $n$, we denote by $TrM$ its trace on the
hyperplane
defined by $x_n$, and by $Res\, M$ the corresponding
residual
scheme, which is again a model of dimension $n$. We
define more generally $Tr^{(p)}M$ and $Res^{(p)}M$
for any nonnegative integer $p$: with the notations introduced above, we
set $Tr^{(p)}M := I_p$,
and define $Res^{(p)}M$ to be the model corresponding to the
ideal
$$I_0\oplus I_1x_n\oplus \cdots \oplus I_{p-1}x_n^{p-1}\oplus
I_{p+1}x_n^{p}\oplus\dots \oplus I_{q+1}x_n^{q}\oplus\dots$$

If $M_1, \dots , M_{\ell}$ are   models, we say that their
disjoint union
${\bf M}$ is a multi-model. If ${\bf p}=(p_1, \dots, p_{\ell})$
is a multi-integer, we define $Tr^{\bf p}{\bf M}$ to be the disjoint union
of the $Tr^{p_i}M_i$ and $Res^{\bf p}{\bf M}$ to be the disjoint union
of the $Res^{p_i}M_i$.

The Hilbert scheme
$Hom({\bf M}, X^0)$ is the product $Hom(M_1, X^0)\times \dots \times
Hom(M_{\ell},X^0)$, thus irreducible
(and smooth). Its generic point represents an embedding whose image in $X$
we denote by ${\bf M}_X$.
We denote by
$Hom({\bf M}, X, H)$ the Hilbert scheme of morphisms $f$ from ${\bf M}$ to
$X^0$
whose restrictions to the components $M_i$ are $H$-morphisms. We
call these morphisms $H$-morphisms from ${\bf M}$ to $X$,
and we call $H$-embeddings
those $H$-morphisms which are embeddings. The scheme
$Hom({\bf M}, X, H)$ is the
product
$Hom(M_1, X, H)\times \dots \times
Hom(M_{\ell},X, H)$, thus irreducible
(and smooth). Its generic point is a $H$-embedding whose image
in $X$ we denote by ${\bf M}_{X,H}$.
We are now ready to state and prove our differential Horace lemma:

\begin{prop}\label{prin.cor} Let $X$ be, as above, a reduced projective
variety of dimension $n$, furnished with a line bundle $\fasm L$,
and $H$ a reduced
irreducible positive Cartier divisor on $X$ not contained in the
singular locus of $X$.
 Let $W\subset X$ be a closed
subscheme of $X$ not containing $H$. We denote by $W\upp$ and
$W'$ the trace and residual of $W$ with respect to $H$.
Let ${\bf M}$ be a multi-model of dimension $n$.
and ${\bf p}$ a multi-integer (of the same length).
Suppose
\begin{enumerate}
\item \makebox[2.5cm][l]{{\bf Dime}} \hspace{5ex} $H^0(H,\fasm
I_{W\upp\cup Tr^{\bf p}{\bf M}_H}\otimes
\fasm L\vert H)=0$
\item \makebox[2.5cm][l]{{\bf Degue}} \hspace{5ex } $H^0(X,\fasm
I_{W\up\cup Res^{\bf p}{\bf M}_{X,H}}\otimes \fasm L(-H))=0.$
\end{enumerate}
Then $H^0(X,\fasm I_{W\cup {\bf M}_X}\otimes \fasm L)=0$.
\end{prop}

\noindent{\bf Proof.} We may suppose that $k$ is algebraically closed.
By semi-continuity, there exist rational points $z_1, \dots, z_{\ell}$ in $H$,
and corresponding $H$-embeddings $e_i: \mbox{Spec}\, k[[x_1, \dots, x_n]]
\rightarrow X$, with $z_i$ as image of the closed point,
allowing us to rewrite our
first assumption as follows:

\makebox[1.5cm][l]{{\bf Dime}} \hspace{5ex} $H^0(H,\fasm
I_{W\upp\cup e_1(Tr^{p_1}{M_1}) \cup \dots \cup
e_{\ell}(Tr^{p_{\ell}}M_{\ell})}\otimes
\fasm L\vert H)=0.$

Similarly, there exist rational points $z'_1, \dots, z'_{\ell}$ in $H$,
and corresponding $H$-embeddings $e'_i: \mbox{Spec}\, k[[x_1, \dots, x_n]]
\rightarrow X$, with $z'_i$ as image of the closed point,
allowing us to rewrite our
second assumption as follows:

\makebox[1.5cm][l]{{\bf Degue}} \hspace{5ex } 
$H^0(X,\fasm
I_{W\up\cup e'_1(Res^{p_1}{M_1}) \cup \dots \cup
e'_{\ell}(Res^{p_{\ell}}M_{\ell})}\otimes \fasm L(-H))=0.$

Since the two conditions do not interfere with one another, we may even
suppose $z_i=z'_i$. Let us now show that it's possible to obtain $e_i=e'_i$.
Using induction on $i$,
all we need to prove is the following statement: given such a point $z$,
a vertically
graded model $M$ in $\mbox{Spec}\, k[[x_1, \dots, x_n]]$, an integer $p$
and two non-empty open
subschemes $S\upp$
in  $Hom(Tr^pM, H, z)$ and $S\up$ in $Hom(Res^pM, X, H, z)$, there exists a
$H$-embedding
$e: \mbox{Spec}\, k[[x_1, \dots, x_n]]
\rightarrow \mbox{Spec}\, \hat {\fasm O}_{X, z} \rightarrow X$ whose
restriction to $Tr^pM$ is in $S\upp$ and whose restriction to
$Res^pM$ is in $S\up$. This statement is easily proven, using the fact
that the restriction maps
$$Hom (M, X, H, z) \rightarrow Hom (Res^p\,M, X, H, z)$$ 
and 
$$Hom (M, X, H, z) \rightarrow Hom (Tr\,M, H, z)
\rightarrow Hom (Tr^p\,M, H, z)$$
are dominant, and the fact
that any element of $Hom (M, X, H, z)$ extends as a $H$-morphism from
$\mbox{Spec}\, k[[x_1, \dots, x_n]]$ to
$ X$.

The proposition then follows by applying our formal lemma \ref{fl} with
$V= H^0(X, \fasm I_W \otimes \fasm L)$. $\hfill \square$

\begin{rem} Our formal lemma can give more accurate information. For instance
if $X$ is a projective space and $H$ a hyperplane, we may handle {\em linear}
embeddings of (multi-)models in a similar way to that used for general
embeddings
of (multi-)models.
\end{rem}

We now make explicit the particular case of the proposition
where the $M_i$ are points (of various multiplicities).
\begin{cor} \label {prin.lem}
 Let $X$ be, as above, a reduced projective
variety of dimension $n$, furnished with a line bundle $\fasm L$,
and $H$ a reduced
irreducible positive Cartier divisor on $X$ not contained in the
singular locus of $X$.
Let $W$ be a closed subscheme of $X$ not containing $H$. Let
$P_1,\ldots ,P_r$ be generic points of $X$, $Q_1, \ldots ,Q_r$
generic points in $H$ and $m_1,\ldots ,m_r$ a sequence of positive integers.
Then
$H^0(\fasm I_{W\cup P_1^{m_1}\cup \cdots \cup P_r^{m_r}}\otimes
\fasm L)=0$ if the following two conditions are satisfied
(see 2.2 for the notation $D_m$):
\be
\item \makebox[2.5cm][l]{{\bf Dime}} \hspace{5ex} $H^0(X,\fasm
I_{W\upp\cup Q_1\cup \cdots \cup Q_r}\, \otimes
\fasm L|_H)=0$
\item \makebox[2.5cm][l]{{\bf Degue}} \hspace{5ex } $H^0(X,\fasm
I_{W\up\cup D_{m_1}(Q_1)\cup \cdots \cup D_{m_r}(Q_r)}\otimes
\fasm L(-H))=0$
\ee
\end{cor}

\noindent{\bf Proof.} This is just proposition 9.1 with
$p_{i}=m_{i}-1$.$\hfill \square$

\begin{rem}
The lemma \ref{diff} is obtained by taking $H=G_a$ in the
previous
corollary.
\end{rem}


%
$$\br{lll}\mbox{J. Alexander}&\quad\quad\quad&
\mbox{A. Hirschowitz}\\
\mbox{Universit\'e d'Angers}&&\mbox{Universit\'e de Nice}\\
\mbox{jea@tonton.univ-angers.fr}&&\mbox{ah@math.unice.fr}\er$$
\end{document}